# A Triple Inequality with Series and Improper Integrals

Florentin Smarandache
Department of Mathematics
University of New Mexico
Gallup, NM 87301, USA

**Abstract**.

As a consequence of the Integral Test we find a triple inequality which bounds up and down both a series with respect to its corresponding improper integral, and reciprocally an improper integral with respect to its corresponding series.

**1. Introduction**.

Before going in details to this triple inequality, we recall the well-known
Integral Test that applies to positive term series:
For all $x \geq 1$ let $f(x)$ be a positive continuous and decreasing function such that $f(n) = a_n$ for $n \geq 1$. Then:

$$\sum_{n=1}^{\infty} a_n \text{ and } \int_1^{\infty} f(x)dx \qquad (1)$$

either both converge or both diverge.

Following the proof of the Integral Test one easily deduces our inequality.

**2. Triple Inequality with Series and Improper Integrals.**

Let's first make the below notations:

$$S = \sum_{n=1}^{\infty} a_n , \qquad (2)$$

$$I = \int_1^{\infty} f(x)dx . \qquad (3)$$

We have the following
Theorem (Triple Inequality with Series and Improper Integrals):

For all x ≥ 1 let f(x) be a positive continuous and decreasing function such that f(n) = $a_n$ for n ≥ 1. Then:

$$S - f(1) \leq I \leq S \leq I + f(1) \tag{4}$$

Proof.

We consider the closed interval [1, n] the function f is defined on split into n-1 unit subintervals [1, 2], [2, 3], …, [n-1, n], and afterwards the total area of the rectangles of width 1 and length f(k), for 2 ≤ k ≤ n, inscribed into the surface generated by the function f and limited by the x-axis and the vertical lines x = 1 and x = n:

$$S_{\inf} = \sum_{k=2}^{n} f(k) = f(2) + f(3) + \ldots + f(n) \quad \text{[inferior sum]} \tag{5}$$

and respectively the total area of the rectangles of width 1 and length f(k), for 1 ≤ k ≤ n-1, inscribed into the surface generated by the function f and limited by the x-axis and the vertical lines x = 1 and x = n:

$$S_{\sup} = \sum_{k=1}^{n-1} f(k) = f(1) + f(3) + \ldots + f(n-1) \quad \text{[superior sum]} \tag{6}$$

But the value of the improper integral $\int_{1}^{\infty} f(x)dx$ is in between these two summations:

$$S_n - f(1) = S_{\inf} \leq \int_{1}^{n} f(x)dx \leq S_{\sup} = S_{n-1} \tag{7}$$

where

$$S_n = \sum_{k=1}^{n} f(k). \tag{8}$$

Now in (7) computing the limit when n ↦ ∞ one gets a double inequality which bounds up and down an improper integral with respect to its corresponding series:

$$S - f(1) \leq I \leq S \tag{9}$$

And from this one has

$$S \leq I + f(1) \tag{10}$$

Therefore, combining (9) and (10) we obtain our triple inequality:

$$S - f(1) \leq I \leq S \leq I + f(1)$$

As a consequence of this, one has a double inequality which bounds up and down a series with respect to its corresponding improper integral, similarly to (9):

$$I \leq S \leq I + f(1) \tag{11}$$

Another approximation will be:

$$S_n \leq S \leq S_n + I_n \tag{12}$$

where

$$I_n = \int_n^\infty f(x)dx \text{ for } n \geq 1 \tag{13}$$

and $I_1 = I$, $S_1 = a_1 = f(1)$.
The bigger is n the more accurate bounding for S.

These inequalities hold even if both the series S and improper integral I are divergent (their values are infinite). According to the Integral Test if one is infinite the other one is also infinite.

**3. An Application.**

Apply the Triple Inequality to bound up and down the series:

$$S = \sum_{k=1}^{\infty} \frac{1}{k^{\wedge}2+4} \tag{14}$$

The function $f(x) = \dfrac{1}{x^{\wedge}2+4}$ is positive continuous and decreasing for $x \geq 1$. Its corresponding improper integral is:

$$I = \int_1^\infty \frac{1}{x^{\wedge}2+4}dx = \lim_{b \to \infty} \int_1^b \frac{1}{x^{\wedge}2+4}dx = \lim_{b \to \infty} [\frac{1}{2} \arctan \frac{x}{2}]_1^b$$

$$= \frac{1}{2} \lim_{b \to \infty} (\arctan \frac{b}{2} - \arctan \frac{1}{2}) = \frac{1}{2}(\frac{\pi}{2} - \arctan 0.5) \approx 0.553574.$$

Hence:

$$0.553574 = I \leq S \leq I + f(1) = 0.553574 + 1/(1^{\wedge}2 + 4) = 0.753574$$

or

$$0.553574 \leq S \leq 0.753574.$$

With a TI-92 calculator we approximate series (14) summing its first 1,000 terms and we get:

$$S_{1000} = \sum(1/(x^2+4),x,1,1000) = 0.659404.$$

Sure the more terms we take the better approach for the series we obtain.

In a similar way one can bound up and down an improper integral with respect to its corresponding series.

**Reference**:

R. Larson, R. P. Hostetler, B. H. Edwards, with assistance of D. E. Heyd, Calculus / Early Transcendental Functions, Houghton Mifflin Co., Boston, New York, 1999.